\documentclass{elsart}%
\usepackage{amsmath}
\usepackage{graphicx}
\usepackage{float}
\usepackage{indentfirst}
\usepackage{verbatim}
\usepackage{algorithm}
\usepackage{algpseudocode}
\usepackage{float}
\usepackage{multicol}
\usepackage{listings}
\usepackage{subfigure}
\usepackage{epsfig}
\usepackage{amsfonts}
\usepackage{amssymb}
\usepackage{epstopdf}%
\everymath{\displaystyle}
\begin{document}
\begin{frontmatter}
\title{Computation of Local Time of Reflecting Brownian Motion and Probabilistic Representation  of the Neumann Problem}
\author[UNCC]{Yijing Zhou},
\author[UNCC]{Wei Cai},
\author[NW]{Elton Hsu}
\address[UNCC]{Department of Mathematics and Statistics,
University of North Carolina at Charlotte, Charlotte, NC 28223-0001}
\address[NW]{Department of Mathematics, Northwestern University,
Evanston, IL 60521}
\bigskip
{\bf Suggested Running Head:}
\\
Local Time of Reflecting Brownian Motion and Probabilistic Representation of the Neumann Problem
\\
\bigskip
{\bf Corresponding Author: }
\\
Prof. Wei Cai \\
Department of Mathematics and Statistics, \\
University of North Carolina at
Charlotte, \\
Charlotte, NC 28223-0001 \\
Phone: 704-687-0628, Fax: 704-687-6415, \\
Email: wcai@uncc.edu
\newpage
\begin{abstract}
In this paper, we propose numerical methods for computing the boundary local time
of reflecting Brownian motion (RBM) in $R^3$ and its use in the probabilistic representation of the solution
of the Laplace equation with the Neumann boundary condition. Approximations of the RBM based on a walk-on-spheres (WOS)
and random walk on lattices are discussed and tested for sampling the RBM
paths and their applicability in finding accurate approximation of the local time and
discretization of the probabilistic formula. Numerical tests for several types of domains
(cube, sphere, and ellipsoid) have shown the convergence of the numerical methods as
the length of the RBM path and number of paths sampled increase.
\end{abstract}
\end{frontmatter}

\section{Introduction}

Traditionally numerical solutions of boundary value problems for partial
differential equations (PDEs) are obtained by using finite difference, finite
element or boundary element methods with both space and/or time
discretizations. This usually requires spatial mesh fine enough to ensure
accuracy, which results in considerable storage space requirement and
computation time. Moreover, the solution process is global, namely, the
solutions of the PDEs have to be found together at all mesh points. However,
in many scientific and engineering applications, local solutions are sometimes
all we need, such as the local electrostatic potential on a molecular surface
where molecular binding activities are most likely to occur or the stress
field at specific locations where the materials are susceptive to failures.
Therefore, it is of practical importance to have a numerical approach which
can give a local solution of the PDEs at some locality of our choice. In the
case of elliptic PDEs, this kind of local numerical method can be constructed
using the well-known probabilistic representation and the associated
Feynman-Kac formula \cite{[16]}\cite{[17]}, which relate It\^o diffusion paths
to the solution of the elliptic PDEs. By sampling the diffusion paths, the
evaluation of the solution at any point in the domain can be done through an
averaging process of the boundary (Dirichlet or Neumann) data under some given
measure on the boundary. Moreover, this method avoids the expensive mesh
generations required by mesh-based methods mentioned above \cite{[15]}.

Our previous work \cite{[4]} , using the Feynman-Kac formula for the Laplace
equation with Dirichlet data, has produced a local method to compute the DtN
(Dirichlet-to-Neumann) mapping for the Laplace operator. In this paper, we
will focus on solving the following Neumann boundary value problem of the
elliptic PDE using a probabilistic approach:
\begin{equation}
\left\{
\begin{aligned} \Delta u &=f, \  on \ D\\ \frac{\partial u}{\partial n}&=\phi, \  on\  \partial D\\ \end{aligned}\right.
, \label{eq1}%
\end{equation}
where $D$ is a bounded domain in $R^{3}$, $\Delta$ is the Laplace operator,
$f$ is a measurable function and $\phi$ is a bounded measurable function on
the boundary $\partial D$ satisfying $\int_{\partial D}\phi d\sigma=\int_{D} f
dx$. Equation $(1)$ becomes the Laplace equation when $f=0$, which is the
subject of our work.

The PDE (\ref{eq1}) originates from either the Poisson equation for
electrostatic potentials \cite{[19]}, an implicit time discretization of the
heat equation or the momentum equation of the Naiver-Stokes equation with an
additional lower term in the latter cases. Historically, Brownian motion (BM)
has been used in solving PDEs due to its effectiveness and easy implementation
regardless of dimensions \cite{[10]}. The well-known probabilistic
representation can solve the elliptic equation with the Dirichlet boundary
condition by using the first exit time $\tau_{D}$ of BM, i.e.,
\begin{equation}
u(x)=E^{x}(\phi(x_{\tau_{D}}))+E^{x}\left[  \int_{0}^{\tau_{D}}f(X_{t}%
)dt\right]  . \label{eq3}%
\end{equation}
In the above formula, only the values at the hitting positions on the boundary
are used in the computation of the mathematical expectation (average) to
obtain $u(x)$. Taking the idea of killed Brownian motion \cite{[4]} in
junction with Monte Carlo methods, we can easily obtain an estimate of $u(x)$.

However, for the Neumann problem to be studied here, in contrast to the
Brownian motion in (\ref{eq3}), reflecting Brownian motion (RBM) will be
needed to produce a similar probabilistic solution to (\ref{eq1}). This theory
has been developed in \cite{[1]} by employing the concept of the boundary
local time whose one dimensional predecessor was introduced by L\'{e}vy in
\cite{[3]}. In \cite{[1]}, the boundary local time of a one dimensional BM was
extended to high dimensions and an explicit form, shown in (\ref{eq13}), was
obtained for domains with smooth boundaries. It should be noted that the
boundary local time is related to the Skorohod equation \cite{[14]} and plays
a significant role in the theoretical development of the probabilistic
approach to the Neumann problem.

One-dimensional local time of Brownian motion has been studied by many authors
\cite{[3]}\cite{[11]}\cite{[13]}\cite{[14]}. For higher dimensions, similar
results have been found by Brosamler \cite{[2]}. Morillon \cite{[8]} gave a
modified Feynman-Kac formula for the Poisson problem with various boundary
conditions, algorithms based on random walk on a grid have been proposed.
However, numerical algorithms for computing local time in $R^{3}$ based on a
rigorous probabilistic theory has not been done in the literature. It is the
objective of this paper to obtain practical numerical methods for computing
the local time of RBM in three dimensions and apply the resulting numerical
methods to implement computationally the probabilistic representation for the
Neumann problem.

The rest of the paper is organized as follows. In section 2, we give some
background information on the Skorohod problem which is the key to the Neumann
problem. In section 3, an explicit probabilistic solution to the Neumann
problem will be given. In section 4, a walk-on-spheres (WOS) method is
reviewed and discussed for its application for RBM. In section 5, a numerical
method, the WOS combined with a Monte Carlo method, is proposed for an
approximation to the Neumann problem. Numerical results for cubic, spherical
and ellipsoidal domains will be given in Section 6. Finally, we draw
conclusions from our Monte Carlo simulations and discuss possible further work.

\section{ Skorohod Problem, RBM and Boundary Local Time}

For the sake of completeness, we first give the definitions of Brownian motion
and reflecting Brownian motion in $R^{d}$.

\begin{defn}
\textrm{Brownian motion}: A Brownian motion $B(t)=(B_{1}(t),B_{2}%
(t),...,B_{d}(t))$ in $R^{d}$ is a set of $d$ independent stochastic processes
with the following properties: for $1\le i\le d$,

\begin{enumerate}
\item (Normal increments) $B_{i}(t)-B_{i}(s)$ has a normal distribution with
mean 0 and variance $t-s$.

\item (Independence of increments) $B_{i}(t)-B_{i}(s)$ is independent of the
past, i.e., of $B_{u}$, $0\leq u\leq s$.

\item (Continuity of paths) $B_{i}(t),t\geq0$ is a continuous function of $t$.
\end{enumerate}
\end{defn}

\bigskip

\begin{defn}
\textrm{Skorohod equation:} Assume $D$ is a bounded domain in $R^{d}$ with a
$C^{2}$ boundary. Let $f(t)$ be a (continuous) path in $R^{d}$ with
$f(0)\in\bar{D}$. A pair $(\xi_{t},L_{t})$ is a solution to the Skorohod
equation $S(f;D)$ if the following conditions are satisfied:

\begin{enumerate}
\item $\xi$ is a path in $\bar{D}$;

\item $L(t)$ is a nondecreasing function which increases only when $\xi
\in\partial D$, namely,
\begin{equation}
L(t)=\int_{0}^{t}I_{\partial D}(\xi(s))L(ds); \label{eq5}%
\end{equation}

\item The Skorohod equation holds:
\begin{equation}
S(f;D):\qquad\ \xi(t)=f(t)-\frac{1}{2}\int_{0}^{t}n(\xi(s))L(ds), \label{eq7}%
\end{equation}
where $n(x)$ stands for the outward unit normal vector at $x\in\partial D$.
\end{enumerate}
\end{defn}

\begin{rem}
\ In Definition 2, the smoothness constraint on $D$ can be relaxed to bounded
domains with $C^{1}$ boundaries, which however will only guarantee the
existence of $(\ref{eq7})$. But for a domain $D$ with a $C^{2}$ boundary, the
solution will be unique. Obviously, $(\xi_{t},L_{t})$ is continuous in the
sense that each component is continuous.
\end{rem}

If $f(t)$ is replaced by the standard Brownian motion (BM) $B_{t}$, the
corresponding $\xi_{t}$ will be a standard reflecting Brownian motion (RBM)
$X_{t}$. Just as the name suggests, a reflecting BM (RBM) behaves like a BM as
long as its path remains inside the domain $D$, but it will be reflected back
inwardly along the normal direction of the boundary when the path attempts to
pass through the boundary. The fact that $X_{t}$ is a diffusion process with
the Neumann boundary condition can be proven by using a martingale formulation
and showing that $X_{t}$ is the solution to the corresponding martingale
problem with the Neumann boundary condition \cite{[1]}. The result gives an
intuitive and direct way to construct RBM from BM. This construction will be
discussed in detail in Section 5.


Next we will review the concept of boundary local time $L(t)$ for a RBM, which
in a sense is a measure of the amount of time a RBM spends near the boundary
and at the same time the frequency that a RBM hits the boundary.
We have the following properties of $L(t)$:

\begin{enumerate}
\item[(a)] It is the unique continuous nondecreasing process that appears in
the Skorohod equation (\ref{eq7}) \cite{[1]};

\item[(b)] It measures the amount of time the standard reflecting Brownian
motion $X_{t}$ spending in a vanishing neighborhood of the boundary within the
period $[0,t]$. If $D$ has a $C^{3}$ boundary, then
\begin{equation}
L(t)\equiv\lim_{\epsilon\rightarrow0}\frac{\int_{0}^{t}I_{D_{\epsilon}}%
(X_{s})ds}{\epsilon}, \label{eq11}%
\end{equation}
where $D_{\epsilon}$ is a strip region of width $\epsilon$ containing
$\partial D$ and $D_{\epsilon}\subset D$. This limit exists both in $L^{2}$
and $P^{x}$-$a.s$. for any $x\in\overline{D}$;

\item[(c)] $L(t)$ is a continuous additive functional (CAF) \cite{[2]} which
satisfies the additivity property \cite{[1]}\cite{[2]}\cite{[13]}\cite{[14]}:
$A_{t+s}=A_{s}+A_{t}(\theta_{s})$. Here $\theta_{s}$ denotes the shift
operator along the paths. The additivity property of $L(t)$ can be seen as follows:

From the definition in (\ref{eq11}), we have
\[
\left\{
\begin{aligned} L(t+s)&= \lim_{\epsilon \rightarrow 0}\frac{\int_{0}^{t+s}I_{D_{\epsilon}}(X_{\tau})d\tau}{\epsilon}\\ L(s)&= \lim_{\epsilon \rightarrow 0}\frac{\int_{0}^{s}I_{D_{\epsilon}}(X_{\tau})d\tau}{\epsilon}\\ L(t)\circ\theta_s&= \lim_{\epsilon \rightarrow 0}\frac{\int_{0}^{t}I_{D_{\epsilon}}(X_{\tau+s})d\tau}{\epsilon}=\lim_{\epsilon \rightarrow 0}\frac{\int_{s}^{t+s}I_{D_{\epsilon}}(X_{\tau})d\tau}{\epsilon} \end{aligned}\right.
,
\]
therefore,
\[
L(t+s)=L(t)+L(t)\circ\theta_{s},
\]
which shows that $L(t)$ satisfies the additivity property.
\end{enumerate}

For one-dimensional case, much existing literature devoted to the study of
local times for Brownian motion and more general diffusion processes. It is
well known that for one dimensional Brownian motion starting from the origin,
the local time $L(t)$ of RBM and $\max_{s\leq t}B(s)$ have the same
distribution as stochastic processes. Hence, valuable properties of RBM can be
drawn by just observing $\max_{s\leq t}B(s)$. But this is not true in general
in higher dimensions. However, we have the following explicit formula for
$L(t)$ derived in \cite{[1]},
\begin{equation}
L(t)=\sqrt{\frac{\pi}{2}}\int_{0}^{t}I_{\partial D}(X_{s})\sqrt{ds},
\label{eq13}%
\end{equation}
where the the right-hand side of (\ref{eq13}) is understood as the limit of
\begin{equation}
\sum_{i=1}^{n-1}%
\smash{\displaystyle\max_{s\in\Delta_i}I_{\partial D}(X_s)\sqrt{|\Delta_i|}},\quad
\smash{\displaystyle\max_i}|\Delta_{i}|\rightarrow0, \label{eq15}%
\end{equation}
where $\Delta=\{\Delta_{i}\}$ is a partition of the interval $[0,t]$ and each
$\Delta_{i}$ is an element in $\Delta$. We will discuss the implementation of
both (\ref{eq11}) and (\ref{eq13}) in Section 5. \newline

\section{Neumann Problem}

We will consider the elliptic PDE in $R^{3}$ with a Neumann boundary
condition
\begin{equation}
\left\{
\begin{aligned} \left(\frac{\Delta}{2}+q\right)u &=0, \  on \ D\\ \frac{\partial u}{\partial n}&=\phi, \  on\  \partial D\\ \end{aligned}\right.
. \label{eq17}%
\end{equation}
When the bottom of the spectrum of the operator $\Delta/2+q$ is negative a
probablistic solution of $(\ref{eq17})$ is given by%
\begin{equation}
u(x)=\frac{1}{2}E^{x}\left[  \int_{0}^{\infty}e_{q}(t)\phi(X_{t})L(dt)\right]
, \label{eq19}%
\end{equation}
where $X_{t}$ is a RBM starting at $x$ and $e_{q}(t)$ is the Feynman-Kac
functional \cite{[1]}
\[
e_{q}(t)=\exp\left[  \int_{0}^{t}q(X_{s})\,ds\right]  .
\]

From the definition of the local time in (\ref{eq11}), we have the following
approximation for small $\epsilon$
\begin{equation}
L(t)\approx\frac{\int_{0}^{t}I_{D_{\epsilon}}(X_{s})ds}{\epsilon}.
\label{eq21}%
\end{equation}
Plugging $(\ref{eq21})$ into $(\ref{eq19})$, we have
\begin{equation}
u(x)\approx\frac{1}{2\epsilon}E^{x}\left[  \int_{0}^{\infty}e_{q}(t)\phi
(X_{t})\int_{t}^{t+dt}I_{D_{\epsilon}}(X_{s})ds\right]  . \label{eq23}%
\end{equation}

The solution defined in $(\ref{eq19})$ should be understood as a weak solution
for the classical PDE $(\ref{eq17})$. The proof of the equivalence of
($\ref{eq19}$)\ with a classical solution is done by using a martingale
formulation \cite{[1]}. If the weak solution satisfies some smoothness
condition \cite{[1]}\cite{[2]}, it can be shown that it is also a classical
solution to the Neumann problem. This formula is the basis for our numerical
approximations to the Neumann problem $(\ref{eq17})$. To compute the
expectation in the formula, we rely on Monte Carlo random samplings to
simulate Brownian paths and then take the average.

In the present work, as we only consider the Laplace equation where $q=0$,
therefore,
\begin{equation}
u(x)\approx\frac{1}{2\epsilon}E^{x}\left[  \int_{0}^{\infty}\phi(X_{t}%
)\int_{s}^{s+dt}I_{D_{\epsilon}}(X_{s})ds\right]  , \label{eq25}%
\end{equation}
and we will show how this formula is implemented with the Monte Carlo and WOS
methods in section 5. \newline

\begin{rem}
Comparing with formula $(\ref{eq19})$, we find that the probabilistic
solutions to the Laplace operator with the Dirichlet boundary condition has a
very similar form (referring to $(\ref{eq29})$). In the Dirichlet case, killed
Brownian paths were sampled by running random walks until the latter are
absorbed on the boundary and $u(x)$ is evaluated as an average of the
Dirichlet values at the first hitting positions on the boundary, namely,
$u(x)=E^{x}\left[  \phi(X_{\tau_{D}})\right]  $ where $\phi$ is the Dirichlet
boundary data. On the other hand, for the Neumann condition, while $u(x)$ is
also given as a weighted average of the Neumann data at hitting positions of
RBM on the boundary, the weight is related to the boundary local time of RBM.
This is a noteworthy point when we compare the probabilistic solutions of the
two boundary value problems and try to understand the formula in
$(\ref{eq19})$.
\end{rem}

\section{Method of Walk on Spheres (WOS)}

Random walk on spheres (WOS) method was first proposed by M\"{u}ller
\cite{[7]}, which can solve the Dirichlet problem for the Laplace operator
efficiently. Here we will first briefly review this method and then show how
it can be adapted for RBM and the Neumann problem.

For a general linear elliptic problem with a Dirichlet boundary condition,
\begin{equation}
\begin{aligned} & L(u) =\sum_{i=1}^{n}b_{i}(x)\frac{\partial u}{\partial x_{i}} + \sum_{i,j=1}^{n}a_{ij}(x)\frac{\partial^{2} u}{\partial x_{i}\partial x_{j}} = f(x), x\in D,\\ & u|_{\partial D} = \phi(x), x\in\partial D.\\ \end{aligned} \label{eq27}%
\end{equation}
The probabilistic representation of the solution is (\cite{[16]}\cite{[17]})
\begin{equation}
u(x)=E^{x}(\phi(x_{\tau_{D}}))+E^{x}\left[  \int_{0}^{\tau_{D}}f(X_{t}%
)dt\right]  , \label{eq29}%
\end{equation}
where $X_{t}(w)$ is an It\^{o} diffusion defined by
\begin{equation}
dX_{t}=b(X_{t})dt+\alpha(X_{t})dB_{t}, \label{eq31}%
\end{equation}
and $B_{t}(w)$ is the Brownian motion, $[a_{ij}]=\frac{1}{2}\alpha
(x)\alpha^{T}(x)$.

The expectation in (\ref{eq29}) is taken over all sample paths starting from
$x$ and $\tau_{D}$ is the first exit time for the domain $D$. This
representation holds true for general linear elliptic PDEs. For the Neumann
boundary condition, similar formulas can be obtained \cite{[8]}. However
different measures on the boundary $\partial D$ will be used in the
mathematical expectation.

In order to illustrate the WOS method for the Dirichlet problem, let us
consider the Laplace equation where $f=0,a_{ij}=\delta_{ij}$ and $b_{i}=0$ in
(\ref{eq27}) and the It\^{o} diffusion is then simply the standard Brownian
motion with no drift. The solution to the Laplace equation can be rewritten in
terms of a measure $\mu_{D}^{x}$ defined on the boundary $\partial D$,
\begin{equation}
u(x)=E^{x}(\phi(X_{\tau_{D}}))=\int_{\partial D}\phi(y)d\mu_{D}^{x},
\label{eq33}%
\end{equation}
where $\mu_{D}^{x}$ is the so-called harmonic measure\ defined by
\begin{equation}
\mu_{D}^{x}(F)=P^{x}\left\{  X_{\tau_{D}}\in F\right\}  ,F\subset\partial
D,x\in D. \label{eq35}%
\end{equation}
It can be shown that the harmonic measure is related to the Green's function
for the domain with a homogeneous boundary condition,
\begin{equation}
\left\{
\begin{aligned} -\Delta g(x,y) &= \delta(x-y), \  &x\in D,\\ g(x,y) &= 0, \  &x\in\partial D\\ \end{aligned}\right.
. \label{eq37}%
\end{equation}
By the third Green's identity,
\begin{equation}
u(x)=\int_{\partial D}\left[  u(y)\frac{\partial g(y,x)}{\partial
n}-g(y,x)\frac{\partial u}{\partial n}(y)\right]  dS_{y}, \label{eq39}%
\end{equation}
and using the zero boundary condition of $g$, we have
\begin{equation}
u(x)=\int_{\partial D}u(y)\frac{\partial g(y,x)}{\partial n}dS_{y}.
\label{eq41}%
\end{equation}
Thus, the hitting probability $\mu_{D}^{x}([y,y+dS_{y}])$ is equivalent to
$p(x,y)dS_{y}$. Comparing (\ref{eq33}) with (\ref{eq41}), we can see that
\begin{equation}
p(\mathbf{x},\mathbf{y})=-\frac{\partial g(x,y)}{\partial n_{y}}. \label{eq43}%
\end{equation}
For instance, the Green's function for a ball for this purpose is given as
\begin{equation}
g(x,y)=-\frac{1}{4\pi|x-y|}+-\frac{1}{4\pi|x-y^{\ast}|}, \label{eq44}%
\end{equation}
where $y^{\ast}$ is the inversion point of $y$ with respect to the sphere
\cite{[4]}.

If the starting point $x$ of a Brownian motion is at the center of a ball, the
probability of the BM exiting a portion of the boundary of the ball will be
proportional to the portion's area. It is known that all sample functions of
Brownian motion processes starting in the domain intersects the boundary
$\partial D$ almost surely \cite{[7]}. Therefore, sampling a Brownian path by
drawing balls within the domain, regardless of how the path navigates in the
interior before hitting the boundary, can significantly reduce the path
sampling time. To be more specific, given a starting point $x$ inside the
domain $D$, we simply draw a ball of largest possible radius fully contained
in $D$ and then the next location of the Brownian path on the surface of the
ball can be sampled, using a uniform distribution on the sphere, say at
$x_{1}$. Treat $x_{1}$ as the new starting point, draw a second ball fully
contained in $D$, make a jump from $x_{1}$ to $x_{2}$ on the surface of the
second ball as before. Repeat this procedure until the path hits a absorption
$\epsilon$-shell of the domain \cite{[5]}. When this happens, we assume that
the path has hit the boundary $\partial D$ (see Figure 1(a) for an illustration).

\begin{figure}[ptb]
{\large \centering   \subfigure[WOS within the domain]{
\label{fig:subfig:a}     \includegraphics[width=0.37\textwidth]{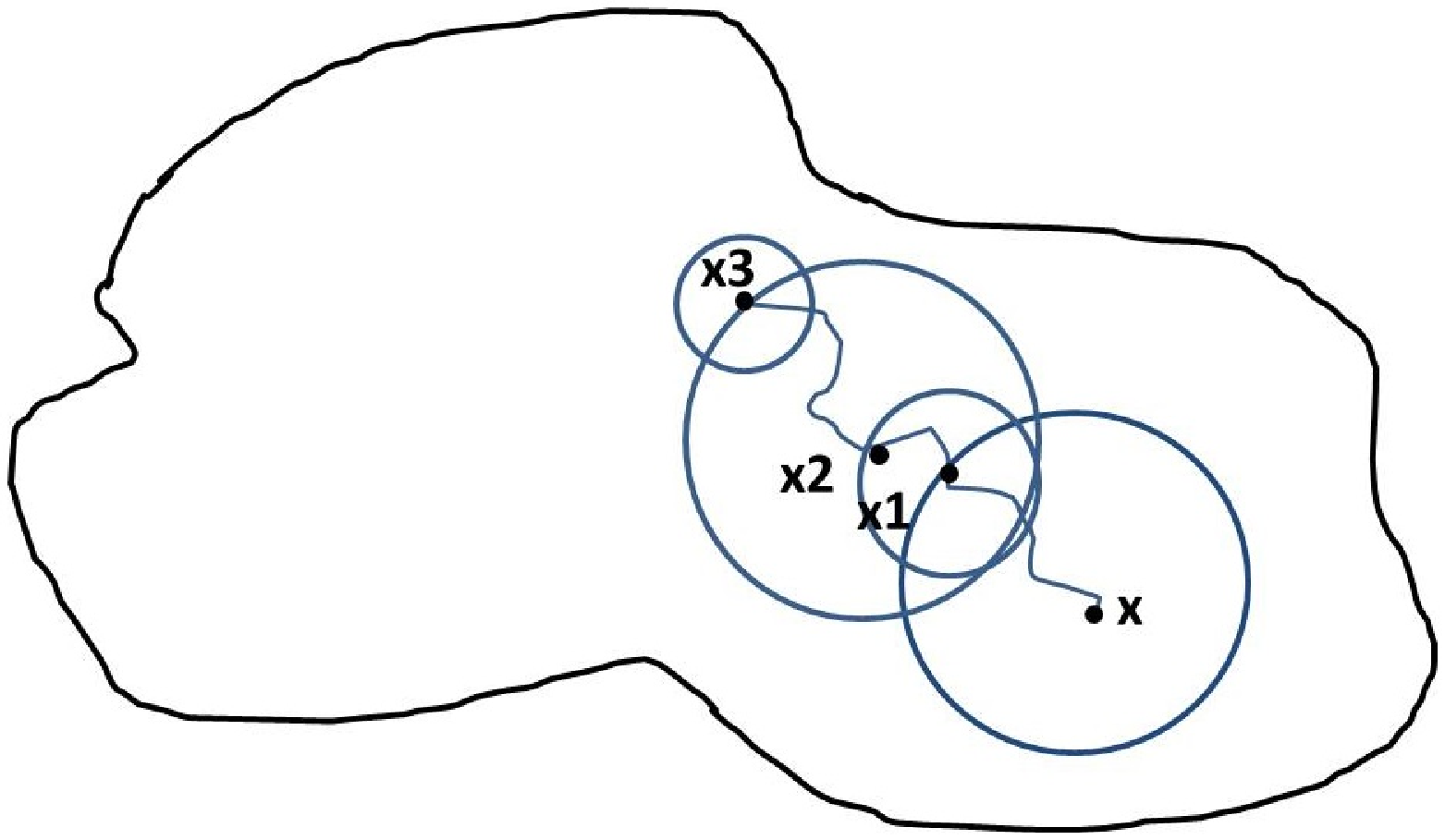}}
\hspace{1in}
\subfigure[WOS (with a maximal step size for each jump) within the domain]{
\label{fig:subfig:b}     \includegraphics[width=0.35\textwidth]{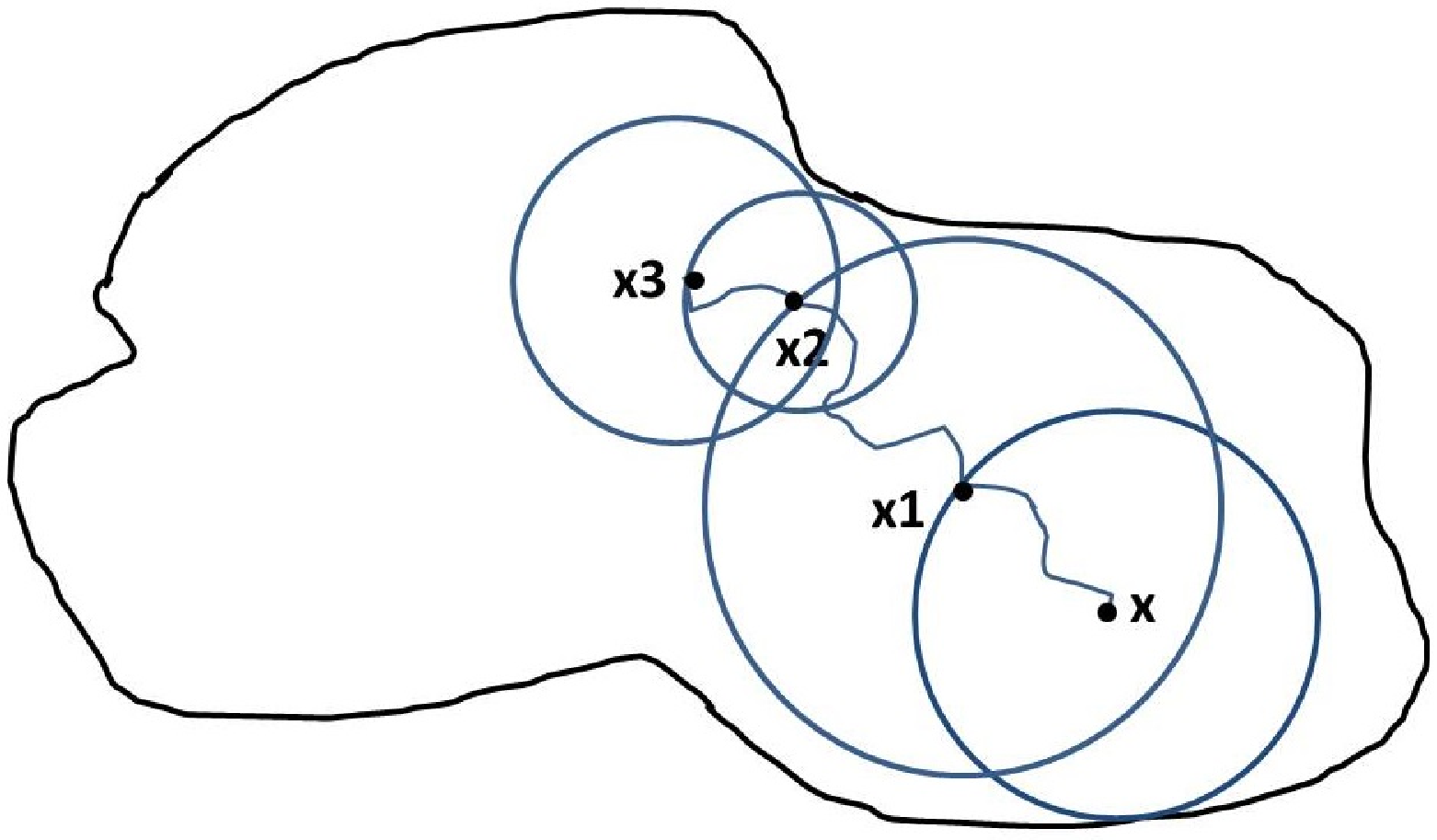}}
}\caption{Walk on Spheres method}%
\label{fig:subfig}%
\end{figure}

Next we define an estimator of (\ref{eq29}) by
\begin{equation}
u(x)\approx\frac{1}{N}\sum_{i=1}^{N}u(x_{i}),
\end{equation}
where $N$ is the number of Brownian paths sampled and $x_{i}$ is the first
hitting point of each path on the boundary. Using a jump size (radius of the
ball) $\delta$ on each step for the WOS, we expect to take $O(1/\delta^{2})$
steps for a Brownian path to reach the boundary \cite{[6]}. To speed up,
maximum possible size for each step would allow faster first hitting on the
boundary. Most of the numerical results in this paper will use the WOS
approach as illustrated in Figure 1(b).

\section{ Numerical Methods}

\subsection{Simulation of reflecting Brownian paths}

A standard reflecting Brownian motion path can be constructed by reflecting a
standard Brownian motion path back into the domain whenever it crosses the
boundary. So in principle, the simulation of RBM is reduced to that of BM.

\begin{figure}[ptb]
\centering {\large \includegraphics[width=0.38\textwidth]{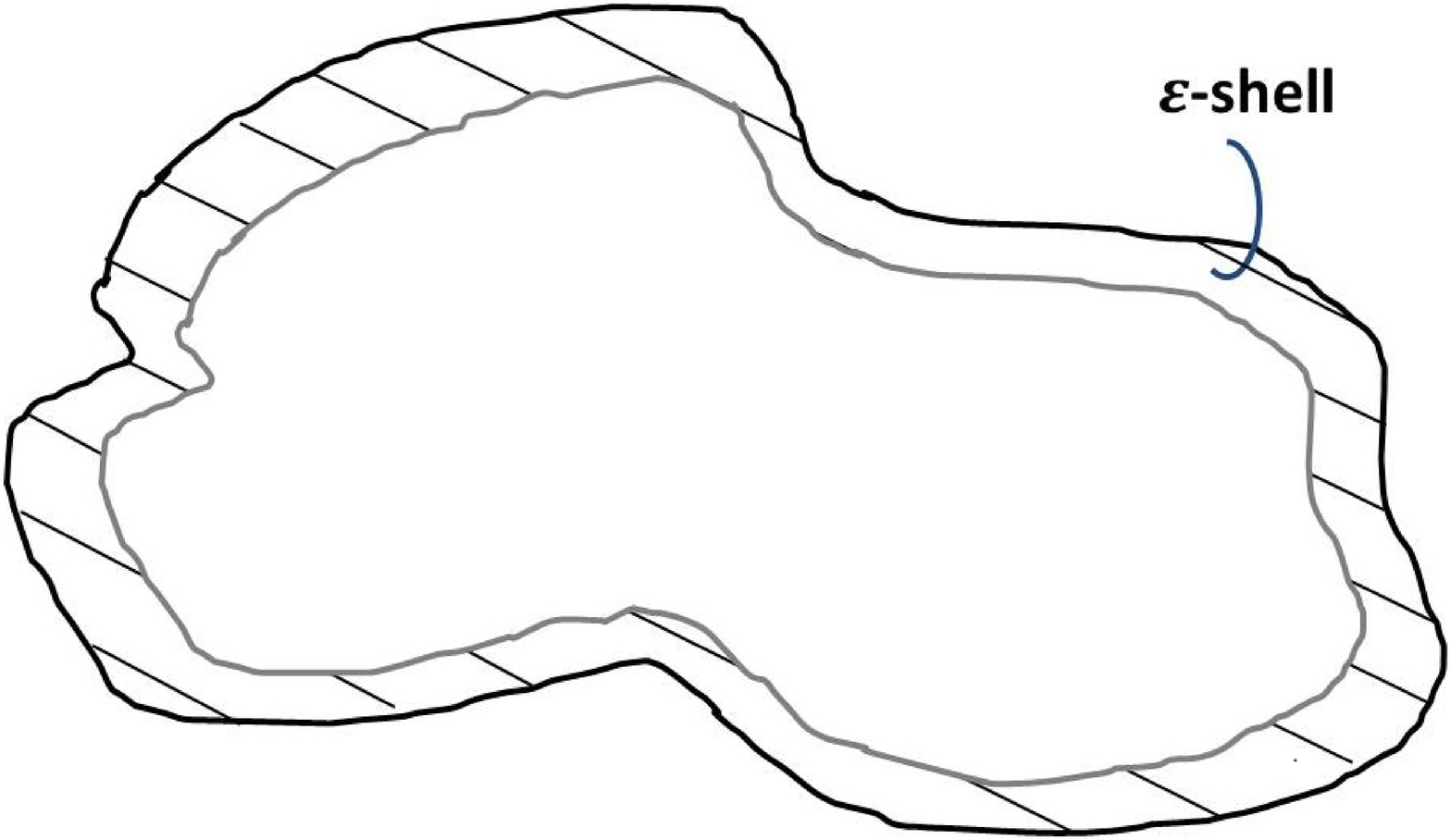}
}\caption{A $\epsilon$-region for a bounded domain in $R^{3}$}%
\end{figure}

It is known that standard Brownian motion can also be constructed as the
scaling limit of a random walk on a lattice so we can model BM by a random
walk with proper scale (see Appendix for details). However, it turns out that
the WOS method is the preferred method to simulate BM for our purpose
\cite{[9]} (see Remark \ref{rm1} for details). As mentioned before, a
$\epsilon$-shell is chosen around the boundary as the termination region in
the Dirichlet case. Here we follow a similar strategy by setting up a
$\epsilon$-region but allowing the process $X_{t}$ to continue moving after
the latter reaches the $\epsilon$-region instead of being absorbed.

Figure 2 shows a strip region with width $\epsilon$ near the boundary is
identified for a bounded domain. In a spherical domain, the $\epsilon$-region
is simply an $\epsilon$-shell near the boundary of width $\epsilon$. Denote
$M_{\epsilon}(D)$ as the $\epsilon$-region and $I(D)$ as the remaining
interior region $D\backslash M_{\epsilon}(D)$.

\begin{figure}[ptb]
\centering {\large \includegraphics[width=0.5\textwidth]{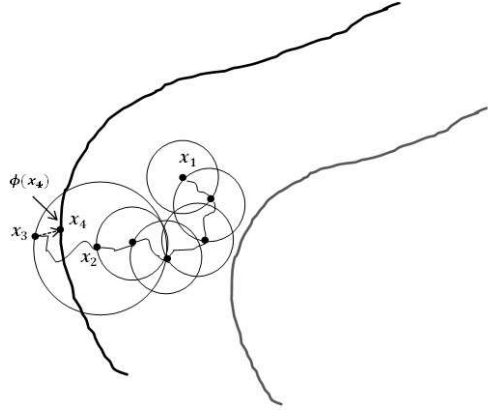}
}\caption{WOS in the $\epsilon$-region. BM path hits $x_{1}$ in $\epsilon
$-region for the first time. Then the radius of sphere is changed to $\Delta
x$, the path continues until it arrives at $x_{2}$ whose distance to $\partial
D$ is smaller or equal to $\Delta x$. Then the radius of the ball is enlarged
to 2$\Delta x$ so that the path has a chance to run out of the domain at
$x_{3}$. If that happens, we pull back $x_{3}$ to $x_{4}$ which is the closest
point to $x_{3}$ on the boundary. Record $\phi(x_{4})$, and continue
WOS-sampling the path starting at $x_{4}$. }%
\end{figure}

Recall the discussion of the WOS in the previous section. For a BM starting at
a point $x$ in the domain, we draw a ball centered at $x$, the Brownian path
will hit the spherical surface with a uniform probability as long as the ball
does not overlap the domain boundary $\partial D$. The balls are constructed
so that the jumps are as large as possible by taking the radius of the ball to
be the distance to the boundary $\partial D$. We repeat this procedure until
the path reaches the region $M_{\epsilon}(D)$. Here, we continue the WOS in
$M_{\epsilon}(D)$ but with a fixed radius $\Delta x$ much smaller than
$\epsilon$. In order to simulate the path of RBM, at some points of the time
the BM path may run out of the domain. For this to happen, the radius of WOS
is increased to $2\Delta x$ when the path is close to boundary at a distance
less than $\Delta x$. In this way, the BM path will have a chance to get out
of the domain, and when that happens, we then pull it back to the nearest
point on the boundary along the normal of the boundary. Afterwards, the BM
path will continue as before.

\begin{figure}[ptb]
\centering {\large \includegraphics[width=0.6\textwidth]{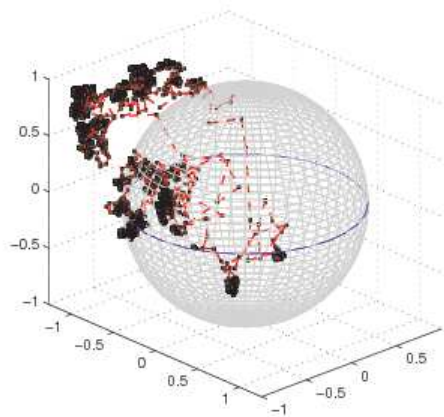}
}\caption{A RBM path with a cube in $R^{3}$}%
\end{figure}

In summary, a reflecting Brownian motion path is simulated by the WOS method
inside $D$. Once it enters the $\epsilon$-region $M_{\epsilon}(D)$, the radius
of WOS changes to a fixed value, either $\Delta x$ or $2\Delta x$, depending
on its current distance of the Brownian particle to the boundary. Once the
path reaches a point on the boundary after the reflection, the radius of WOS
changes back to $\Delta x$. Figure 3 illustrates the movement of RBM in the
$\epsilon$-region $M_{\epsilon}(D)$. As time progresses, we expect the path
hits the boundary at some time instances and lies in either $I(D)$ or
$M_{\epsilon}(D)$ at others. A RBM path is shown in Figure 4 within a cube of
size 2.

\subsection{Computing the boundary local time $L(t)$}

Two equivalent forms of the local time have been given in (\ref{eq11}) and
(\ref{eq13}). Here we will show how the $\epsilon$-region for the construction
of the RBM in Fig. 3 can also be used for the calculation of the local time.
When the $\epsilon$-region is thin enough, i.e. $\epsilon\ll1$, an
approximation of (\ref{eq11}) is given in (\ref{eq21}), which is the
occupation time that RBM $X_{s}$ sojourns within the $\epsilon$-region during
the time interval $[0,t]$. A close look at (\ref{eq21}) reveals that only the
time spent near the boundary is involved and the specific moment when the path
enters the $\epsilon$-region has no effect on the calculation of $L(t)$.

Suppose $x\in D$ is the starting point of a Brownian path, which is simulated
by the WOS method. Once the path enters the $\epsilon$-region, the radius of
WOS is changed to $\Delta x$ or $2\Delta x$. It is known that the elapsed time
$\Delta t$ for a step of a random walk on average is proportional to the
square of the step size, in fact, $\Delta t=(\Delta x)^{2}/d,d=3$ when $\Delta
x$ is small (see Appendix), which also applies to WOS moves (See Remark 6 for
details). Therefore, we can obtain an approximation of the local time $L(dt)$
by counting the number of steps the path spent inside $M_{\epsilon}(D)$
multiplied by the time elapsed for each step, i.e.%
\begin{equation}
L(dt)=L(t_{j}-t_{j-1})\approx\frac{\int_{t_{j-1}}^{t_{j}}I_{D_{\epsilon}%
}(X_{s})ds}{\epsilon}=(n_{t_{j}}-n_{t_{j-1}})\frac{(\Delta x)^{2}}{3\epsilon},
\label{eq45}%
\end{equation}
where $n_{t_{j}}-n_{t_{j-1}}$ is the number of steps that WOS steps remain in
the $\epsilon$-region during the time interval [$t_{j-1},t_{j}].$ Figure 5
gives a sample path of the simulated local time associated with the RBM in
Figure 4.

\begin{rem}
(Alternative way to compute local time $L(t)$ ) From $(\ref{eq13})$, the local
time increases if and only if the RBM path hits the boundary, which implies
that the time before the path hits the boundary makes no contribution to the
increment of the local time. Thus, a WOS method with a changing radius can
also be used with $(\ref{eq13})$. Specifically, we divide the time interval
$[0,t]$ into to $N$ small subintervals of equal length. In each $[t_{i}%
,t_{i+1}]$ the Brownian path will move $2\Delta x$ or $\Delta x$ with the WOS
method when the current path lies within a distance less or more than $\Delta
x$ to the boundary. If the path hits or crosses the boundary within
$[t_{i},t_{i+1}]$, then $L(t)$ will increase by $\sqrt{\pi/2}\sqrt
{t_{i+1}-t_{i}}$.
\end{rem}

\begin{figure}[ptb]
\centering {\large \includegraphics[width=0.6\textwidth]{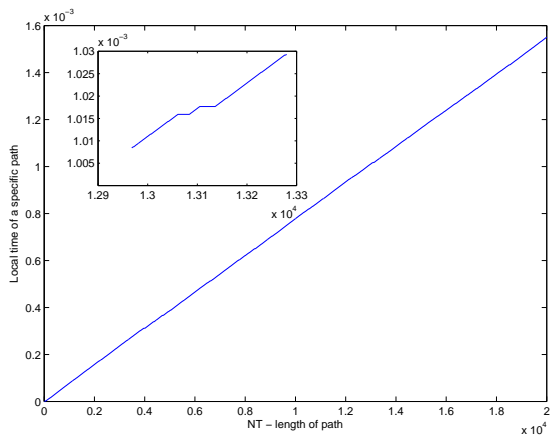}
}\caption{Boundary local time (\ref{eq45}) increases when the path runs into
the region $M_{\epsilon}(D)$. The insert shows the piecewise linear profile of
the local time path with flat level regions. The path of $L(t)$ is a
nondecreasing function and increases at a rate lower than $O(NT)$. $NT$ is the
length of the path.}%
\end{figure}

\begin{figure}[ptb]
\centering {\large \includegraphics[width=0.5\textwidth]{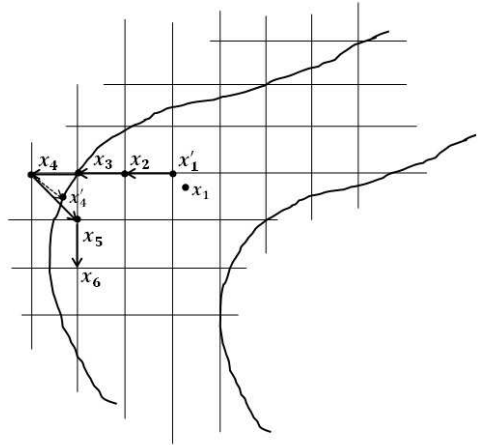}
}\caption{Random walks on the $\epsilon$-region. A BM path hits $x_{1}\in
M_{\epsilon}(D)$ by the WOS method. Replace $x_{1}$ by the nearest grid point
$x_{1}^{\prime}$. Then several steps of random walks will make a path as
$x_{2}\rightarrow x_{3}\rightarrow x_{4}$. Since $x_{4}\notin D$, we push it
back along the normal line (dash arrow) to $x_{4}^{\prime}$ then replace it by
the closest grid point within domain (solid arrow) $x_{5}$. Here path crosses
the boundary at $x_{4}^{\prime}\in\partial D$. Then continue the random walk
as usual at $x_{6}$.}%
\end{figure}

\begin{rem}
\label{rm1} (Approximating RBM by WOS or random walks on a lattice - a
comparison) There are two ways to find approximation to Brownian paths inside
the region $M_{\epsilon}(D)$ and construct their reflections once they get out
of the boundary. One way is by using the WOS approach as described in Section
5.1. The other is in fact to use a random walk on a lattice inside
$M_{\epsilon}(D).$ In the second approach, as illustrated in Fig. 6, a grid
mesh is set up over $M_{\epsilon}(D)$ and the random walk takes a one-step
walk on the lattice until the path goes out of the domain and then it will be
pushed back to the nearest lattice point inside $M_{\epsilon}(D)$. And the
elapsed time for a $\Delta x$ walk is on average $(\Delta x)^{2}/3$ as shown
in the Appendix. The boundary local time $L(t)$ can be still calculated as in
$(\ref{eq45})$. The problem with this approach is that a Brownian motion
actually should have equal probability to go in all directions in the space
while a random walk on the lattice only considers six directions in $R^{3}$.
This limitation was found in our numerical tests to lead to insufficient
accuracy in simulating reflecting Brownian motions for our purpose.

Meanwhile, the WOS method in the $\epsilon$-region $M_{\epsilon}(D)$ has a
fixed radius $\Delta x$, which enables us to calculate the boundary local time
by $(\ref{eq45})$ since the elapsed time of a $\Delta x$ move in $R^{3}$ on
average still remains to be $(\Delta x)^{2}/3$. \ This conclusion can be
heuristically justified by considering points on the sphere are the linear
combination of the directions along the three axes, which implies that the
average time that the path hits the sphere with radius $\Delta x$ should also
be the same. As discussed before, if the path comes within a distance very
close to the boundary, say less than $\Delta x$, the radius of the WOS method
is increased to $2\Delta x$ so that it will have a chance to run out of the
domain and then be pushed back to the nearest point on the boundary to affect
a hit of the RBM on the boundary.
\end{rem}

\subsection{Probabilistic representation for the Neumann problem}

Finally, with the boundary local time of RBM available, we can come to the
approximation of the Neumann problem solution $u(x)$ using the probabilistic
approach (\ref{eq25}). First of all, we will need to truncate the infinite
time duration required for the RBM path $X_{t}$ in (\ref{eq25}) to a finite
extent for computer simulations. The exact length of truncation will have to
be numerically determined by increasing the length until a convergence is
confirmed (namely, the approximation to $u(x)$ does not improve within a
prescribed error tolerance between two different choices of truncation times
under same number of sampled paths). Assume that the time period is limited to
from $0$ to $T$, then by a Monte Carlo sampling of the RBM paths, an
approximation of (\ref{eq25}) will be
\begin{equation}
\widetilde{u}(x)=\frac{1}{2\epsilon}\sum_{i=1}^{N}\left[  \int_{0}^{T}%
\phi(X_{t}^{i})I_{\partial D}(X_{t}^{i})\int_{t}^{t+dt}I_{D_{\epsilon}}%
(X_{s}^{i})ds\right]  , \label{eq47}%
\end{equation}
where $X_{t}^{i},i=1,...,N$ are stochastic processes sampled according to the
law of RBM.

Next, let us see how the RBM can be incorporated into the representation
formula once its path is obtained.

Associate the time interval $[0,T]$ with the number of steps $NT$ of a
sampling path, $NT$ will give the total length of each path. Then, the
integral inside the square bracket in (\ref{eq47}) can be transformed into
$j$
\begin{equation}
\sum_{j^{\prime}=1}^{NT}\left(  \phi(X_{t_{j}}^{i})I_{\partial D}(X_{t_{j}%
}^{i})\int_{t_{j-1}}^{t_{j}}I_{D_{\epsilon}}(X_{s}^{i})ds\right)  ,
\label{eq49}%
\end{equation}
where $j^{\prime}$ stands for the $j^{\prime}-$th step the WOS method has
taken, and $j$ indicates a step for which $X_{t_{j}}^{i}\in\partial D$.

As the integral in (\ref{eq49}) is in fact the occupation time as shown in
(\ref{eq45}), (\ref{eq49}) becomes
\begin{equation}
\sum_{j^{\prime}=1}^{NT}\left(  \phi(X_{t_{j}}^{i})I_{\partial D}(X_{t_{j}%
}^{i})(n_{t_{j}}-n_{t_{j-1}})\frac{(\Delta x)^{2}}{3}\right)  . \label{eq51}%
\end{equation}

As a result, an approximation to the PDE solution $\widetilde{u}(x)$ becomes
\begin{equation}
\widetilde{u}(x)=\frac{1}{2\epsilon}\sum_{i=1}^{N}\left[  \sum_{j^{\prime}%
=1}^{NT}\left(  \phi(X_{t_{j}}^{i})I_{\partial D}(X_{t_{j}}^{i})(n_{t_{j}%
}-n_{t_{j-1}})\frac{(\Delta x)^{2}}{3}\right)  \right]  . \label{eq53}%
\end{equation}

Theoretically speaking, $\epsilon$ should be chosen much larger than $\Delta
x$. Here, we take $\epsilon=k\Delta x$, $k$ $>1$ is an integer, which will
increase as $\Delta x$ vanishes to zero. Then, (\ref{eq53}) reduces to
\begin{equation}
\begin{aligned} \widetilde{u}(x) &= \frac{1}{2k\Delta x} \sum_{i=1}^{N}\left[\sum_{j=1}^{NT}\left(\phi(X_{t_j}^i)I_{\partial D}(X_{t_j}^i) (n_{t_j}-n_{t_{j-1}})\frac{(\Delta x)^2}{3}\right)\right]\\ &=\frac{\Delta x}{6k}\sum_{i=1}^{N}\left[\sum_{j=1}^{NT}\left(\phi(X_{t_j}^i)I_{\partial D}(X_{t_j}^i) (n_{t_j}-n_{t_{j-1}})\right)\right],\\ \end{aligned} \label{eq55}%
\end{equation}
which is the final numerical algorithm for the Neumann problem. In the
following we present the general implementation of this numerical algorithm.

Let $x$ be any interior point in $D$ where the solution $u(x)$ for the Neumann
problem is sought. First, we define the $\epsilon$-region $M_{\epsilon}(D)$
near the boundary. For each one of $N$ RBM paths, the following procedure will
be executed until the length of the path reaches a prescribed length given by
$NT\cdot\Delta x$:

\begin{enumerate}
\item If $x\notin M_{\epsilon}(D)$, predict next point of the path by the WOS
with a maximum possible radius until the path locates near the boundary within
a certain given distance $\epsilon,$ say $\epsilon=5\Delta x$ (hit the
$\epsilon$-region $M_{\epsilon}(D)$). If $x\in M_{\epsilon}(D)$, $l(t_{i})=1$;
otherwise, $l(t_{i})=0$. Here $l(t)$ is the unit increment of $L(t)$ at time
$t$.

\item If $x\in M_{\epsilon}(D)$, use the WOS method with a fixed radius
$\Delta x$ to predict the next location for Brownian path. Then, execute one
of the two options:
\end{enumerate}

\qquad\emph{Option 1}. If the path happens to hit the domain boundary
$\partial D$ at $x_{t_{i}}$, record $\phi(x_{t_{i}})$.

\qquad\emph{Option 2}. If the path passes crosses the domain boundary
$\partial D$, then pull the path back along the normal to the nearest point on
the boundary. Record the Neumann value at the boundary location.

Due to the independence of the paths simulated with the Monte Carlo method, we
can run a large number of paths simultaneously on a computer with many cores
in a perfectly parallel manner, and then collect all the data at the end of
the simulation to compute the average. \textbf{Algorithm 1} gives a
pseudo-code for the numerical realization of implementing the WOS in both
$I(D)$ and $M_{\epsilon}(D)$ regions. \newline

\begin{algorithm}[H]
\begin{algorithmic}
\State \textbf{Data:} Select integers $N$ and $NT$, a starting point $X_0\in D$, step size $h$
and $\epsilon$-region $M_\epsilon(D)$ near the boundary.
\State \textbf{Output} An approximation of $u(X_0)$.\\
\State \textbf{Initialization} $L[NT], v[NT], u[N]$, $X=X_0$, $i\leftarrow1$ and $j\leftarrow1$;\\
\State \textbf{While} {$i\le N$} \textbf{do}
\State \hspace*{2ex}   Set $S_i=0$.
\State \hspace*{2ex}   \textbf{While} $j\le NT$ \textbf{do}
\State \hspace*{4ex}\quad   \textbf{If} $X\in I(D)$ \textbf{then} \quad /* If the path has not touched the $\epsilon$-region */
\State \hspace*{4ex}\qquad Set $L[j]\leftarrow$ 0; \quad /*Increment of local time at each step. */
\State \hspace*{4ex}\qquad Set $r\leftarrow d(X,\partial D)$;\quad  /* Find the distance to the boundary */
\State \hspace*{4ex}\qquad Randomly choose a point $X_1$ on $B(X,r)$ then set $X\leftarrow X_1$.
\State \hspace*{4ex}\quad \textbf{Else}\quad  /* The path enters the $\epsilon$-region  */
\State \hspace*{4ex}\qquad $L[j]\leftarrow1$; \quad /*local time increases */
\State \hspace*{4ex}\qquad Set $r\leftarrow$ h (2h);\quad   /* If $d(X,\partial D)> h$ or =0 ($0< d(X,\partial D)\le h$) */
\State \hspace*{4ex}\qquad Randomly choose a point $X_1$ on $B(X,r)$ then set $X\leftarrow X_1$.
\State \hspace*{4ex}\qquad \textbf{If} $X\notin \bar{D}$, \textbf{then}
\State \hspace*{6ex}\qquad  Find $X_j$ to be the nearest point on $\partial D$ to X and pull $X$ back
\State \hspace*{6ex}\qquad onto $\partial D$  at $X_j$;
\State \hspace*{6ex}\qquad Set $X\leftarrow X_j$;
\State \hspace*{6ex}\qquad Set $v[j]\leftarrow\phi(X_j)$
\State \hspace*{4ex}\qquad \textbf{End}
\State \hspace*{4ex}\quad \textbf{End}
\State \hspace*{4ex}\quad  $j\leftarrow j+1$;
\State \hspace*{2ex}\textbf{End}
\State \hspace*{2ex} count $\leftarrow$ 0;
\State \hspace*{2ex} \textbf{For} k=1:NT
\State \hspace*{4ex} count $\leftarrow$ count + $L[k]$;
\State \hspace*{4ex} \textbf{If} $v[k]\sim=0$ \textbf{then}
\State \hspace*{4ex} $u[i]\leftarrow u[i]+\phi(X_k)\cdot$count;
\State \hspace*{4ex} count $\leftarrow$ 0;
\State \hspace*{2ex} \textbf{End}
\State \hspace*{2ex}$i\leftarrow i+1$;
\State  \textbf{End}
\State \textbf{Return} $\widetilde{u}(X_0)=h\sum_{k=1}^{N}u[k]/N/(6k)$
\end{algorithmic}
\end{algorithm}

\begin{center}
\vspace*{-15pt} Algorithm 1: The algorithm for the probabilistic solution of
the Laplace equation with the Neumann boundary condition
\end{center}

\section{Numerical results}

In this section, we give the numerical results for the Neumann problem in
cubic, spherical and ellipsoid domains.

To monitor the accuracy of the numerical approximation of the solutions, we
select a circle inside the domain, where the solution of the PDE $u(x)$ will
be found by the proposed numerical methods, defined by
\begin{equation}
\{(x,y,z)^{T}=(r\cos\theta_{1}\sin\theta_{2},r\sin\theta_{1}\sin\theta
_{2},r\cos\theta_{2})^{T}\} \label{eq57}%
\end{equation}
with $r=0.6$, $\theta_{1}=0:k\cdot2\pi/30:2\pi$, $\theta_{2}=\pi/4$ with
$k=1,...,15$. In addition, a line segment will also be selected as the
locations to monitor the numerical solution, the endpoints of the segment are
$(0.4,0.4,0.6)^{T}$ and $(0.1,0,0)^{T}$, respectively. Fifteen uniformly
spaced points on the circle or the line are chosen as the locations for
computing the numerical solutions.

The true solution of the Neumann problem (\ref{eq1}) with the corresponding
Neumann boundary data is
\begin{equation}
u(x)=\sin3x\sin4y\ e^{5z}+5. \label{eq58}%
\end{equation}

In the figures of numerical results given below, the blue curves are the true
solutions and the red-circle ones are the approximations. The numerical
solutions are shifted by a constant so they agree with the exact solution at
one point as the Neumann problem is only unique up to an arbitrary additive
constant. ``Err''indicates the relative error of the approximations.

\subsection{Cube domain and test on the length of the path}

A cube domain of size 2 is selected to test the choice of the number of paths
and the length of the paths (truncation time duration $T$) in the numerical
formula (\ref{eq55}).

The step-size $\Delta x=0.0005$ is used as the radius of the WOS inside the
$\epsilon$-region $M_{\epsilon}(D)$, namely, the step-size of the random walk
approximation of the RBM near the boundary. The number of paths is taken as
$N=2e5.$ Two choices for the path length parameter $NP=2.7e4$ and $NP=3e4$ are
compared to gauge the convergence of the numerical formula (\ref{eq55}) in
terms of the path truncation. Figures 7 and 8 shows the solution and the
relative errors in both cases, which indicates that $NP=3e4$ will be
sufficient to give an error around 5\% as shown in Fig. 8.

In the rest of the numerical tests, we will set the number of path $N=2e5,$
and number of steps for each path $NP=3e4$.

\begin{figure}[ptb]
{\large \centering   \subfigure[$\epsilon$ = 6$\Delta x$, Err = 10.71\%]{
\label{fig:subfig:a}
\includegraphics[width=0.38\textwidth]{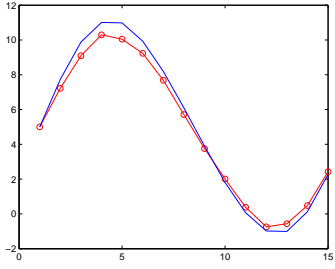}} \hspace{1in}
\subfigure[$\epsilon$ = 7$\Delta x$, Err = 12.19\%]{
\label{fig:subfig:b}
\includegraphics[width=0.38\textwidth]{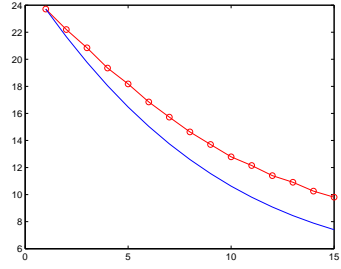}} }\caption{Cubic domain:
number of paths $N=2e5$, and number of steps for each path $NP=2.7e4$. (Left)
Solution on the circle defined in (\ref{eq58}), (right) solution on a line
segment.}%
\label{fig:subfig}%
\end{figure}

\begin{figure}[ptb]
{\large \centering   \subfigure[$\epsilon$ = 6$\Delta x$, Err = 5.35\%]{
\label{fig:subfig:a}
\includegraphics[width=0.38\textwidth]{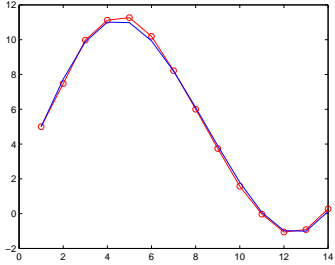}} \hspace{1in}
\subfigure[$\epsilon$ = 7$\Delta x$, Err = 5.85\%]{
\label{fig:subfig:b}
\includegraphics[width=0.38\textwidth]{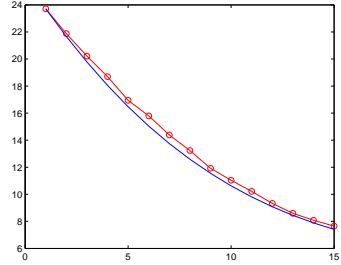}} }\caption{Cubic domain:
number of paths $N=2e5,$ and number of steps for each path $NP=3e4$. (Left)
Solution on the circle defined in (\ref{eq58}), (right) solution on a line
segment.}%
\label{fig:subfig}%
\end{figure}

\subsection{Spherical domain}

The unit ball is centered at the origin. We set $\Delta x=0.0005$ and adjust
$\epsilon$, similar numerical results are obtained as in the case of the cube
domain. Here, the reflected points of Brownian path are the intersection of
the normal and the domain. Though Figure 8(b) shows some oscillations in the
middle, the overall approximation are within a relative error around $5.85\%$.

\begin{figure}[ptb]
{\large \centering   \subfigure[$\epsilon$ = 5$\Delta x$, Err = 5.13\%]{
\label{fig:subfig:a}
\includegraphics[width=0.38\textwidth]{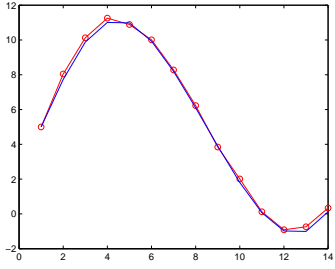}} \hspace{1in}
\subfigure[$\epsilon$ = 5$\Delta x$, Err = 4.03\%]{
\label{fig:subfig:b}
\includegraphics[width=0.38\textwidth]{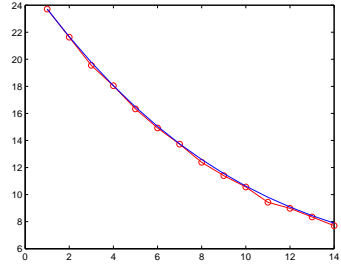}} }\caption{Spherical domain:
number of paths $N=2e5$, and number of steps for each path $NP=3e4$. (Left)
Solution on the circle defined in (\ref{eq58}), (right) solution on a line
segment.}%
\label{fig:subfig}%
\end{figure}

\subsection{Ellipsoid domain}

The ellipse with axis lengths $(3,\ 2,\ 1)$ is centered at the origin. We set
$\Delta x=0.0004$. The numerical results along the circle behave better than
those along the line segment, especially along the tail section of the latter
(Figure 10(b)), which lie closer to the origin $\mathbf{0}$.

\begin{figure}[ptb]
{\large \centering   \subfigure[$\epsilon$ = 5$\Delta x$, Err = 5.75\%]{
\label{fig:subfig:a}
\includegraphics[width=0.38\textwidth]{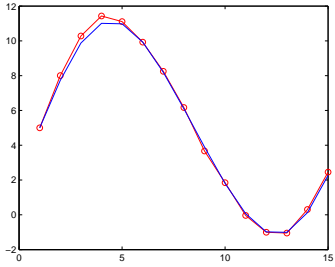}} \hspace{1in}
\subfigure[$\epsilon$ = 5$\Delta x$, Err = 5.12\%]{
\label{fig:subfig:b}
\includegraphics[width=0.38\textwidth]{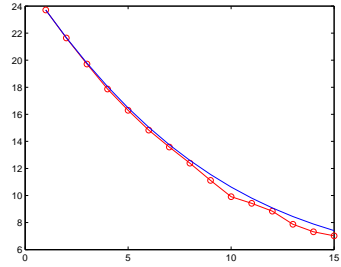}} }\caption{Ellipsoid domain:
number of paths $N=2e5$, and number of steps for each path $NP=3e4$. Left:
solution on the circle defined in (\ref{eq58}). Right: solution on a line
segment.}%
\label{fig:subfig}%
\end{figure}



\section{Conclusions and discussion}

In this paper we have proposed numerical methods for computing the local time
of reflecting Brownian motion and the probabilistic solution of the Laplace
equation with the Neumann boundary condition. Without knowing the complete
trajectories of RBM in space, we are able to use the WOS to sample the RBM and
calculate its local time, based on which a discrete probabilistic
representation (\ref{eq55}) was obtained to produce satisfactory
approximations to the solution of the Neumann problem at one single point.
Numerical results validated the stability and accuracy of the proposed
numerical methods.

In addition, random walk on a lattice was also investigated as an alternative
way to sample RBM. However, numerical experiments show that the numerical
results are inferior to those obtained by the WOS method. A possible reason is
that formula (\ref{eq11}) for the local time is valid for a smooth path while
a random walk approximation of the the Brownian path contains inherent errors.

The local time can also be computed by a mathematically equivalent formula
(\ref{eq13}), for which the implementation is discussed briefly in section
5.2. Again the numerical results based on (\ref{eq13}) are inferior to those
obtained using the original limiting process of L\'{e}vy in \cite{[3]} . This
fact we believe may result from the time discretization error of Brownian
paths especially when long time truncation is employed in the probabilistic representation.

Various issues affecting the accuracy of the proposed numerical methods remain
to be further investigated, such as the number of random walk or WOS steps and
the truncation of duration time $T$ for the paths, the choice of the thickness
for the $\epsilon$-region, the size of $\Delta x$ for the lattice, etc. In
theory, the larger the truncation time $T$, the more accurate is the
probabilistic formula for the Neumann solution. However, for a fixed spatial
mesh size $\Delta x,$ long time integration will result in the accumulation of
time discretization error for the Brownian pathes, thus leading to the
degeneracy of the numerical solutions as verified by our numerical
experiments. Meanwhile, R\'{e}v\'{e}sz \cite{[18]} have proposed some
approximations of local time by other stochastic processes in the case of a
half line. We conjecture such results may still hold in higher dimensions and
progress in this direction will shed light on how to improve the numerical
procedures proposed in this paper.

\section{Appendix \newline}

If the random walk on a lattice as in Fig. 11 is to converge to a continuous
BM, a relationship between $\Delta t$ and $\Delta x$ in $R^{3}$ will be needed
and is shown to be
\begin{equation}
\Delta t=\frac{(\Delta x)^{2}}{3}.
\end{equation}

The following is a proof for this result (See \cite{[12]} for reference). The
density function of standard BM satisfies the following PDE \cite{[1]}
\begin{equation}
\frac{\partial p}{\partial t}=\frac{1}{2}\Delta_{x}p(t,x,y)\ . \label{eq59}%
\end{equation}
By using a central difference scheme and changing $p$ to $v$, equation
(\ref{eq59}) becomes
\begin{equation}
\frac{v_{i,j,k}^{n+1}-v_{i,j,k}^{n}}{\Delta t}=\frac{1}{2}\frac{v_{i+1,j,k}%
^{n}+v_{i-1,j,k}^{n}+v_{i,j+1,k}^{n}+v_{i,j-1,k}^{n}+v_{i,j,k+1}%
^{n}+v_{i,j,k-1}^{n}-6v_{i,j,k}^{n}}{(\Delta x)^{2}}\text{ } . \label{eq61}%
\end{equation}

Reorganizing and letting $\lambda=\Delta t/(2(\Delta x)^{2})$ , we have
\begin{equation}
v_{i,j,k}^{n+1}=\lambda v_{i+1,j,k}^{n}+\lambda v_{i-1,j,k}^{n}+\lambda
v_{i,j+1,k}^{n}+\lambda v_{i,j-1,k}^{n}+\lambda v_{i,j,k+1}^{n}+\lambda
v_{i,j,k-1}^{n}+(1-6\lambda)v_{i,j,k}^{n}\text{ },
\end{equation}
By setting $\lambda=\frac{1}{6}$, we have
\begin{equation}
v_{i,j,k}^{n+1}=\frac{1}{6}v_{i+1,j,k}^{n}+\frac{1}{6}v_{i-1,j,k}^{n}+\frac
{1}{6}v_{i,j+1,k}^{n}+\frac{1}{6}v_{i,j-1,k}^{n}+\frac{1}{6}v_{i,j,k+1}%
^{n}+\frac{1}{6}v_{i,j,k-1}^{n} . \label{eq65}%
\end{equation}

\begin{figure}[ptb]
\centering {\large \includegraphics[width=0.5\textwidth]{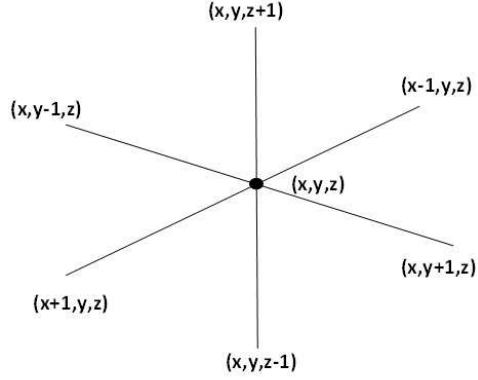}
}\caption{Central difference scheme in $R^{3}$}%
\end{figure}

For the initial condition $\phi$, we have
\begin{equation}
v_{i,j,k}^{n+1}=\sum_{i^{\prime},j^{\prime},k^{\prime}}C_{i^{\prime}%
,j^{\prime},k^{\prime}}\phi\left(  \sum_{l=1}^{n}\overset{\rightarrow}%
{\eta_{l}}\right)  \label{eq67}%
\end{equation}
where
\begin{equation}
\overset{\rightarrow}{\eta_{l}}=\left\{
\begin{aligned} &(-h,0,0)^T, \  &prob = \frac{1}{6}\\ &(h,0,0)^T, \  &prob = \frac{1}{6}\\ &(0,h,0)^T, \  &prob = \frac{1}{6}\\ &(0,-h,0)^T, \  &prob = \frac{1}{6}\\ &(0,0,h)^T, \  &prob = \frac{1}{6}\\ &(0,0,-h)^T, \  &prob = \frac{1}{6}\\ \end{aligned}\right.
, \label{eq69}%
\end{equation}
and
\begin{equation}
\sum_{l=1}^{n}\overset{\rightarrow}{\eta_{l}}=\left(
\begin{aligned} &-n+2i'+i\\ &-n+2j'+j\\ &-n+2k'+k\\ \end{aligned}\right)  h .
\label{eq71}%
\end{equation}
Let $\overset{\rightarrow}{\eta_{l}}=(x_{l},y_{l},z_{l})^{T}$, then
\begin{equation}
x_{l}=\left\{
\begin{aligned} -h, \  &prob=\frac{1}{6}\\ h, \  &prob=\frac{1}{6}\\ 0, \  &prob=\frac{2}{3}\\ \end{aligned}\right.
, \label{eq73}%
\end{equation}
for each $l$. We known that $y_{l}$, $z_{l}$ have the same distribution as
$x_{l}$.

Notice that the covariance between any two of $x_{l}$, $y_{l}$, $z_{l}$ is
zero, i.e. $E(x_{l}y_{l})=0$, $E(y_{l}z_{l})=0$ and $E(x_{l}z_{l})=0$. So
$E(\sum_{i=1}^{n}x_{l}\sum_{i=1}^{n}y_{l})=0$, $E(\sum_{i=1}^{n}y_{l}%
\sum_{i=1}^{n}z_{l})=0$ and $E(\sum_{i=1}^{n}x_{l}\sum_{i=1}^{n}z_{l})=0$.
According to the central limit theorem, we have
\begin{equation}
\sum_{i=1}^{n}x_{l}\overset{D}{=}N\left(  0,\frac{nh^{2}}{3}\right)
\ as\ \ n\rightarrow\infty. \label{eq75}%
\end{equation}
The same assertion holds for $\sum_{i=1}^{n}y_{l}$ and $\sum_{i=1}^{n}z_{l}$.

Since $\lambda=\frac{\Delta t}{2(\Delta x)^{2}}=\frac{1}{6}$, then $h^{2}=3k$
and hence $\frac{nh^{2}}{3}=nk=t$. Therefore $\sum_{i=1}^{n}x_{l}\ \sim
N(0,t)$ as $n\rightarrow\infty$. So are $\sum_{i=1}^{n}y_{l}$ and $\sum
_{i=1}^{n}z_{l}$.

Recall that the covariance between any pair of $\sum_{i=1}^{n}x_{l}$,
$\sum_{i=1}^{n}y_{l}$, and $\sum_{i=1}^{n}z_{l}$ is zero, that $\sum_{i=1}%
^{n}x_{l}$,$\sum_{i=1}^{n}y_{l}$ and $\sum_{i=1}^{n}z_{l}$ are independent
normal random variables. Hence,
\begin{equation}
C_{i^{\prime},j^{\prime},k^{\prime},n}=P\left\{  \sum_{l=1}^{n}\overset
{\rightarrow}{\eta_{l}}=\left(
\begin{aligned} &-n+2i'+i\\ &-n+2j'+j\\ &-n+2k'+k\\ \end{aligned}\right)
h=\left(
\begin{aligned} \sum_{i=1}^{n}x_l\\ \sum_{i=1}^{n}y_l\\ \sum_{i=1}^{n}z_l\\ \end{aligned}\right)
\right\}  \overset{D}{\rightarrow}\frac{1}{(2\pi t)^{3/2}}e^{\frac
{-\Vert\overset{\rightarrow}{x}-\overset{\rightarrow}{x_{0}}\Vert^{2}}{2t}},
\label{eq77}%
\end{equation}
and
\begin{equation}
v_{i,j,k}^{n+1}=\sum_{i^{\prime},j^{\prime},k^{\prime}}C_{i^{\prime}%
,j^{\prime},k^{\prime},n}\phi(\sum_{l=1}^{n}\overset{\rightarrow}{\eta_{l}%
})\rightarrow\iiint_{R^{3}}\frac{1}{(2\pi t)^{3/2}}e^{\frac{-\Vert
\overset{\rightarrow}{x}-\overset{\rightarrow}{x_{0}}\Vert^{2}}{2t}}%
\phi(\overset{\rightarrow}{x})d\overset{\rightarrow}{x}, \label{eq79}%
\end{equation}
which coincides with the density function of the 3-$d$ standard BM.

In conclusion, when $\frac{\Delta t}{2(\Delta x)^{2}}=\frac{1}{6}$, i.e.
$\Delta t=\frac{(\Delta x)^{2}}{3}$ or $\sqrt{dt}=\frac{dx}{\sqrt{3}}$, the
central difference scheme converges to the standard BM in 3-$d$. Generally,
the result can be extended to $d$-dimensional Euclidean space and the result
will be $\Delta t=\frac{(\Delta x)^{2}}{d}$.

\section*{Acknowledgement\newline}

The authors Y.J.Z and W.C. acknowledge the support of the National Science
Foundation (DMS-1315128) and the National Natural Science Foundation of China
(No. 91330110) for the work in this paper. For the present work, the
author E. H. was supported in part by a Simons Foundation Collaboration Grant
for Mathematicians and by a research grant administered through the University
of Science and Technology of China.


\end{document}